\documentclass[a4paper,reqno,10pt]{amsart}

\usepackage{geometry,times}
\usepackage{amsmath,amssymb}

\newtheorem{Theorem}{Theorem}[section]

\newtheorem{Definition}[Theorem]{Definition}
\newtheorem{Corollary}[Theorem]{Corollary}
\newtheorem{Lemma}[Theorem]{Lemma}

\newtheorem{Algorithm}[Theorem]{Algorithm}
\newenvironment{Proof}
{\begin{trivlist}\item[]{{\sc Proof.}}}{\hfill{$\square$}\noindent\end{trivlist}}

\begin{document}
\title{Caps in $\mathbf{\mathbb{Z}_n^2}$}

\author{Sascha Kurz}

\address{Sascha Kurz\\Department of Mathematics, Physics and Informatics\\University of Bayreuth\\Germany}
\email{sascha.kurz@uni-bayreuth.de}

\begin{abstract}
  We consider point sets in $\mathbb{Z}_n^2$ where no three points are on a line -- also called caps or arcs. For the
  determination of caps with maximum cardinality and complete caps with minimum cardinality we provide integer
  linear programming formulations and identify some values for small $n$.
\end{abstract}


\maketitle

\section{Introduction}

\noindent
A $k$-cap in $\mathbb{F}_q^d$ is a subset $A\subseteq \mathbb{F}_q^d$ of size $k$, where no three points are collinear.
A cap is complete if it is not contained in any larger cap. From the combinatorial point of view there are two very natural questions: What is the maximum or minimum size of a cap in $\mathbb{F}_q^d$? If we consider caps in projective spaces over $\mathbb{F}_q$ instead of affine spaces then caps are related to linear codes of minimum distance $4$, see e.~g.{} \cite{1031.51007}.

\medskip

\noindent
For the maximum size $m_2\!\left(\mathbb{F}_q^d\right)$ of a cap in $\mathbb{F}_q^d$ the following values are known \cite{1005.51003,1027.51015,1023.51007}:
\begin{itemize}
 \item[(1)] $m_2\!\left(\mathbb{F}_q^2\right)=\begin{cases}q+2&\text{for even }q,\\q+1&\text{for odd }q.\end{cases}$
 \item[(2)] $m_2\!\left(\mathbb{F}_q^3\right)=q^2$ for $q>2$.
 \item[(3)] $m_2\!\left(\mathbb{F}_2^d\right)=2^d$.
 \item[(4)] $m_2\!\left(\mathbb{F}_3^4\right)=20$, $m_2\!\left(\mathbb{F}_3^5\right)=45$, and
            $m_2\!\left(\mathbb{F}_4^4\right)=40$.
\end{itemize}

\medskip

\noindent
By $n_2\!\left(\mathbb{F}_q^d\right)$ we denote the minimum size of a complete cap in $\mathbb{F}_q^d$. For the projective case in \cite{0987.51007} constructive upper bounds are given. In \cite{1083.51004} the authors consider permutations which they interpret as point sets over $\mathbb{Z}_n^2$ and ask for the minimum number of collinear triples in such configurations. Closely related we ask for the maximum cardinality $\sigma\!\left(\mathbb{Z}_n^2\right)$ of a cap in $\mathbb{Z}_n^2$ where every translation of the two axes contain at most one point.

\medskip

\noindent
In this article we consider similar questions in $\mathbb{Z}_n^2$ instead of $\mathbb{F}_q^d$. For dimension $2$ some authors use the term arcs instead of caps. For $p$ being the smallest prime divisor of an integer $n$ the bounds
$$
  \max\left\{4,\sqrt{2p}+\frac{1}{2}\right\}\le n_2\!\left(\mathbb{Z}_n^2\right)\le\max\left\{4,p+1\right\}
$$
were proven in \cite{pre05050776}. 
For coprime integers $n,m>1$ the bound
$$
  m_2\!\left(\mathbb{Z}_{nm}^2\right)\le\min\left\{n\cdot m_2\!\left(\mathbb{Z}_m^2\right),
  m_2\!\left(\mathbb{Z}_n^2\right)\cdot m\right\}
$$
can be proven.
Whenever the value of $n$ is clear from the context we use the abbreviation $\overline{a}:=a+\mathbb{Z}n$ for integers $a$. By $\sigma\!\left(\mathbb{Z}_n^2\right)$ we denote the maximum
cardinality of a cap in $\mathbb{Z}_n^2$, where each horizontal line $\left(\overline{1},\overline{0}\right)\cdot\mathbb{Z}_n$ and each vertical line $\left(\overline{0},\overline{1}\right)\cdot\mathbb{Z}_n$ contains at most one point. The last conditions model permutations in some sense, see e.~g.{} \cite{1083.51004}.

\subsection{Related work}
The \textit{original} ``no-three-in-line'' problem, introduced by H.~Dudeney in 1917, asks if it is possible to select $2n$ points on the $n$-by-$n$ grid so that no three points are collinear. Currently it is known that for $\varepsilon>0$ and sufficiently large $n$ at least $\left(\frac{3}{2}-\varepsilon\right)$ points can be chosen so that no triple is collinear. Guy conjectures that $\frac{\pi}{\sqrt{3}}n\approx 1.814\, n$ is asymptotically the best possible. In \cite{1118.68106} an analogous question is treated in three-dimensional space. The question for the minimum size of a complete cap in projective planes over finite fields was originally posed by B.~Segre in the late 1950s. In a more general context some authors consider caps (or arcs) over so called projective Hjelmslev planes, see e.~g.{} \cite{1062.51009,0994.51006,1069.51003,pre05202817,dipl_michael,gm2,gm1,1054.51005}. Here we remark that $\mathbb{Z}_{p^r}^2$ is the affine part for the chain ring $\mathbb{Z}_{p^r}$.

\subsection{Our contribution}
In this article we develop an algorithm which can decide whether three given points in $\mathbb{Z}_n^2$ are collinear or not in $O(n\log n)$, given the prime factorization of $n$. We model the problem of the exact determination of $m_2\!\left(\mathbb{Z}_n^2\right)$, $n_2\!\left(\mathbb{Z}_n^2\right)$, and $\sigma\!\left(\mathbb{Z}_n^2\right)$ as integer linear programs. Finally we perform some computer calculations to determine some so far unknown values, e.~g.{} we validate $m_2\!\left(\mathbb{Z}_{25}^2\right)=20$.

\subsection{Organization of the paper}
In Section \ref{sec_collinear} we specify when we consider three points of $\mathbb{Z}_n^2$ to be collinear and develop a fast algorithm which can check collinearity. To this end some cumbersome and technical but elementary calculations have to be executed. In Section \ref{sec_ilp_formulations} we give integer linear programming formulations and in Section \ref{sec_bounds_exact_values} we combine them with some symmetry breaking techniques to determine exact values of $m_2\!\left(\mathbb{Z}_n^2\right)$, $n_2\!\left(\mathbb{Z}_n^2\right)$, and $\sigma\!\left(\mathbb{Z}_n^2\right)$ for small $n$.

\section{Points on a line}
\label{sec_collinear}

\begin{figure}[htp]
  \begin{center}
    \setlength{\unitlength}{0.20cm}
    \begin{picture}(12.2,12.2)
      \multiput(0.15,0)(0,1){13}{\line(1,0){12}}
      \multiput(0.15,0)(1,0){13}{\line(0,1){12}}
      \put(3.65,7.45){\circle*{0.6}}   
      \put(4.65,0.45){\circle*{0.6}}   
      \put(5.65,5.45){\circle*{0.6}}   
      \put(6.65,10.45){\circle*{0.6}}  
      \put(7.65,3.45){\circle*{0.6}}   
      \put(8.65,8.45){\circle*{0.6}}   
      \put(9.65,1.45){\circle*{0.6}}   
      \put(10.65,6.45){\circle*{0.6}}  
      \put(11.65,11.45){\circle*{0.6}} 
      \put(0.65,4.45){\circle*{0.6}}   
      \put(1.65,9.45){\circle*{0.6}}   
      \put(2.65,2.45){\circle*{0.6}}   
    \end{picture}
  \end{center}
  \caption{The line $\left(\overline{3},\overline{7}\right)+\left(\overline{1},\overline{5}\right)\cdot\mathbb{Z}_{12}$
  over $\mathbb{Z}_{12}^2$.}
  \label{fig:example_line}
\end{figure}
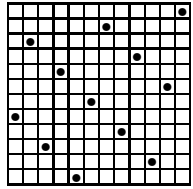

\noindent
A \textbf{line} in $\mathbb{Z}_n^2$ is a translate of a cyclic subgroup of order $n$. We remark that every cyclic subgroup of $\mathbb{Z}_n^2$  is contained in some subgroup of order $n$, see also \cite{pre05050776}. An example is given by the line $\left(\overline{3},\overline{7}\right)+\left(\overline{1},\overline{5}\right)\cdot\mathbb{Z}_{12}$ in $\mathbb{Z}_{12}^2$, see Figure \ref{fig:example_line}. A point $p$ is called \textbf{incident} with a line $l$ if $p\in l$. With this we could define $r$ points to be \textbf{collinear} if they are incident with a common line.

\begin{Lemma}
  \label{lemma_collinear}
  A set of $r$ points $\left(u_i,v_i\right)\in\mathbb{Z}_n^2$ is collinear if there exist
  $a,b,t_1,t_2,w_i\in \mathbb{Z}_n$ with
  $$
    a+w_it_1=u_i\,\,\quad\text{and}\,\,\quad b+w_it_2=v_i
  $$
  for $1\le i\le r$.
\end{Lemma}
\begin{Proof}
  We can write a line $L$ in $\mathbb{Z}_n^2$ as a set
  $
    \Big\{\left(a+wt_1,b+wt_2\right)\mid w\in\mathbb{Z}_n\Big\},
  $
  where $a,b,t_1,t_2$ are arbitrary elements of $\mathbb{Z}_n$.
\end{Proof}

We would like to remark that one could also define a line as the set of solutions $(x,y)\in\mathbb{Z}_n^2$ of where $ax+by=c$, where $a,b,c\in\mathbb{Z}_n$ with $\gcd\left(\widehat{a},\widehat{b},n\right)=1$, see \cite{pre05050776}.

If $n=p$ is a prime then every two points of $\mathbb{Z}_p^2$ uniquely determine a line. This does not hold in general for arbitrary $n$. If $n=p^r$ is a prime power then every two points uniquely determine $p^k$ lines containing those points, where $k$ depends on the points, see e.~g.{} \cite{0994.51006,1054.51005}. It is possible to define neighborhood relations $\sim_k$ for points $a$, $b$ by requiring that there exist at least $p^k$ lines containing $a$ and $b$. For arbitrary $n$ two distinct points are incident with at least one line.

If $n=p$ is a prime then $\mathbb{Z}_n$ is a field and there is a well-known test to check whether three points are collinear or not, which runs in time $O(1)$:

\begin{Lemma}
  \label{lemma_characterization_collinear}
  For a prime $n$ three points $\left(u_1,v_1\right),\left(u_2,v_2\right),\left(u_3,v_3\right)\in\mathbb{Z}_n^2$
  are collinear if and only if
  \begin{eqnarray}
    \label{eqn_collinear}
    \left|\begin{array}{rrr}u_1&v_1&\overline{1}\\u_2&v_2&\overline{1}\\u_3&v_3&\overline{1}\\\end{array}\right|
    =\overline{0}.
  \end{eqnarray}
\end{Lemma}

We remark that in $\mathbb{Z}_8^2$ the points $\left(\overline{0},\overline{0}\right)$, $\left(\overline{2},\overline{4}\right)$, $\left(\overline{4},\overline{4}\right)$ fulfill Equation (\ref{eqn_collinear}) but are not collinear. So in general equation (\ref{eqn_collinear}) is necessary but not sufficient for three points to be collinear.

\begin{Lemma}
  \label{lemma:crit}
  If three points $\left(u_1,v_1\right),\left(u_2,v_2\right),\left(u_3,v_3\right)\in\mathbb{Z}_n^2$ are collinear,
  then
  \begin{eqnarray*}
    \left|\begin{array}{rrr}u_1&v_1&\overline{1}\\u_2&v_2&\overline{1}\\u_3&v_3&\overline{1}\\\end{array}\right|
    =\overline{0}
  \end{eqnarray*}
  holds.
\end{Lemma}
\begin{Proof}
  Due to Lemma \ref{lemma_collinear} we have
  \[
    \left|\begin{array}{rrr}u_1&v_1&\overline{1}\\u_2&v_2&\overline{1}\\u_3&v_3&\overline{1}\\\end{array}\right|=
    \left|\begin{array}{rrr}a+w_1t_1&b+w_1t_2&\overline{1}\\a+w_2t_1&b+w_2t_2&\overline{1}\\a+w_3t_1&b+w_3t_2&
    \overline{1}\\\end{array}\right|=\left|\begin{array}{rrr}w_1t_1&w_1t_2&\overline{1}\\w_2t_1&w_2t_2&\overline{1}\\
    w_3t_1&w_3t_2&\overline{1}\\\end{array}\right|=\overline{0}.
  \]
\end{Proof}

In the following we will develop an algorithm which can decide if three given points in $\mathbb{Z}_n^2$ are collinear in time $O\left(\sum\limits_{i=1}^l\mu_i\right)$ using the prime factorization of $n=\prod\limits_{i=1}^lp_i^{\mu_i}$ as input. We would like to remark that also the Smith normal form can be utilized to obtain an $O(\log n)$-algorithm.

At first we remark that w.l.o.g.{} we can assume that one of three points equals $\left(\overline{0},\overline{0}\right)$:
\begin{Lemma}
  \label{lemma:fix_zero}
  Three points $\left(u_1,v_1\right),\left(u_2,v_2\right),\left(u_3,v_3\right)\in\mathbb{Z}_n^2$ are collinear,
  if and only if there exist $t_1,t_2,w_2,w_3\in \mathbb{Z}_n$ with
  \begin{eqnarray*}
    && w_2t_1=u_2-u_1,\\
    && w_3t_1=u_3-u_1,\\
    && w_2t_2=v_2-v_1,\text{ and}\\
    && w_3t_2=v_3-v_1.
  \end{eqnarray*}
\end{Lemma}
\begin{Proof}
  Since $\left(u_1,v_1\right)$, $\left(u_2,v_2\right)$, $\left(u_3,v_3\right)$ are collinear if and only if
  $\left(0,0\right)$, $\left(u_2-u_1,v_2-v_1\right)$, $\left(u_3-u_1,v_3-v_1\right)$ are collinear, there exist
  $a',b',t_1',t_2',w_1',w_2',w_3'\in \mathbb{Z}_n$ with
  \begin{eqnarray*}
    && a'+w_1't_1'=\overline{0},\\
    && a'+w_2't_1'=u_2-u_1,\\
    && a'+w_3't_1'=u_3-u_1,\\
    && b'+w_1't_2'=\overline{0},\\
    && b'+w_2't_2'=v_2-v_1,\text{ and}\\
    && b'+w_3't_2'=v_3-v_1.
  \end{eqnarray*}

  If we have a solution $t_1,t_2,w_2,w_3\in \mathbb{Z}_n$ of the first equation system then $a'=b'=w_1'=0$, $t_1'=t_1$,
  $t_2'=t_2$, $w_2'=w_2$, $w_3'=w_3$ is a solution of the second equation system.

  If we have a solution $a',b',t_1',t_2',w_1',w_2',w_3'\in \mathbb{Z}_n$ of the second equation system then
  $a'=-w_1't_1'$ and $b'=-w_1't_2'$ holds. With this $t_1=t_1'$, $t_2=t_2'$, $w_2=w_2'-w_1'$,
  $w_3=w_3'-w_1'$ is a solution of the first equation system.
\end{Proof}

Due to the Chinese remainder theorem we have:
\begin{Lemma}
  \label{lemma_decompose}
  If $n=a\cdot b$ with coprime $a$ and $b$, then three points $p_1$, $p_2$, $p_3$ $\in\mathbb{Z}_n^2$ are
  collinear if and only if both projections into $\mathbb{Z}_a^2$ and $\mathbb{Z}_b^2$ give a collinear point set.
\end{Lemma}

Thus it suffices to consider the case where $n=p^r$ is a prime power. For $a\in\mathbb{Z}_n$ let $\hat{a}\in\{1,\dots,n\}$ denote the integer fulfilling $\widehat{a}+\mathbb{Z}_n=a$. If none of the values $u_2,u_3,v_2,v_3$ is invertible in $\mathbb{Z}_{p^r}$ we can apply the following reduction:
\begin{Lemma}
  \label{lemma_reduction}
  If for $r\ge 2$ the prime $p$ divides $\widehat{u}_2$, $\widehat{u}_3$, $\widehat{v}_2$, and $\widehat{v}_3$
  then the points
  $\left(\overline{0},\overline{0}\right)$, $\left(u_2,v_2\right)$ and $\left(u_3,v_3\right)$ are
  collinear in $\mathbb{Z}_{p^r}^2$ if and only if the points $$\left(\overline{0},\overline{0}\right),\,
  \left(\frac{\widehat{u}_2}{p}+\mathbb{Z}p^{r-1},\frac{\widehat{v}_2}{p}+\mathbb{Z}p^{r-1}\right),\text{ and }
  \left(\frac{\widehat{u}_3}{p}+\mathbb{Z}p^{r-1},\frac{\widehat{v}_3}{p}+\mathbb{Z}p^{r-1}\right)$$ are
  collinear in $\mathbb{Z}_{p^{r-1}}^2$.
\end{Lemma}
\begin{Proof}
  We have the following equivalence: If and only if the points $\left(\overline{0},\overline{0}\right)$,
  $\left(u_2,v_2\right)$ and $\left(u_3,v_3\right)$ are collinear in $\mathbb{Z}_{p^r}^2$ then
  due to Lemma \ref{lemma:fix_zero} there exist integers
  $\widetilde{t}_1,\widetilde{t}_2,\widetilde{w}_2,\widetilde{w}_3\in\{1,\dots,p^r\}$ fulfilling
  $$
    p^r\,|\,\widehat{u}_2-\widetilde{w}_2\widetilde{t}_1,\quad
    p^r\,|\,\widehat{u}_3-\widetilde{w}_3\widetilde{t}_1,\quad
    p^r\,|\,\widehat{v}_2-\widetilde{w}_2\widetilde{t}_2,\text{ and }
    p^r\,|\,\widehat{v}_3-\widetilde{w}_3\widetilde{t}_2.
  $$
  Analogously we conclude: If and only if the points $$\left(\overline{0},\overline{0}\right),\,
  \left(\frac{\widehat{u}_2}{p}+\mathbb{Z}p^{r-1},\frac{\widehat{v}_2}{p}+\mathbb{Z}p^{r-1}\right),\text{ and }
  \left(\frac{\widehat{u}_3}{p}+\mathbb{Z}p^{r-1},\frac{\widehat{v}_3}{p}+\mathbb{Z}p^{r-1}\right)$$ are
  collinear in $\mathbb{Z}_{p^{r-1}}^2$ then there exist integers
  $\dot{t}_1,\dot{t}_2,\dot{w}_2,\dot{w}_3\in\{1,\dots,p^{r-1}\}$ fulfilling
  $$
    p^{r-1}\,\Big|\,\frac{\widehat{u}_2}{p}-\dot{w}_2\dot{t}_1,\quad
    p^{r-1}\,\Big|\,\frac{\widehat{u}_3}{p}-\dot{w}_3\dot{t}_1,\quad
    p^{r-1}\,\Big|\,\frac{\widehat{v}_2}{p}-\dot{w}_2\dot{t}_2,\text{ and }
    p^{r-1}\,\Big|\,\frac{\widehat{v}_3}{p}-\dot{w}_3\dot{t}_2.
  $$

  If the tuple $\dot{t}_1,\dot{t}_2,\dot{w}_2,\dot{w}_3\in\{1,\dots,p^{r-1}\}$ is a solution of the second
  system, then $\widetilde{t}_1=\dot{t}_1\cdot p$, $\widetilde{t}_2=\dot{t}_2\cdot p$,
  $\widetilde{w}_2=\dot{w}_2$, $\widetilde{w}_3=\dot{w}_3$ is a solution with
  $\widetilde{t}_1,\widetilde{t}_2,\widetilde{w}_2,\widetilde{w}_3\in\{1,\dots,p^r\}$ of the first system.

  If the tuple $\widetilde{t}_1,\widetilde{t}_2,\widetilde{w}_2,\widetilde{w}_3\in\{1,\dots,p^r\}$ is a solution
  of the first system then $\dot{t}_1=\frac{\widetilde{t}_1}{p}$, $\dot{t}_2=\frac{\widetilde{t}_2}{p}$,
  $\dot{w}_2=\widetilde{w}_2$, $\dot{w}_3=\widetilde{w}_3$ is a solution with
  $\dot{t}_1,\dot{t}_2,\dot{w}_2,\dot{w}_3\in\{1,\dots,p^{r-1}\}$ of the second system.
\end{Proof}

Thus in the following we can confine ourselves to the case where $r\ge 2$ (for $r=1$ we can utilize Lemma \ref{lemma_characterization_collinear}) and at least one of $u_2$, $u_3$, $v_2$, or $v_3$ is invertible in $\mathbb{Z}_{p^r}$. We claim that in such a situation the criterion of Lemma \ref{lemma:crit} is also sufficient for $\left(\overline{0},\overline{0}\right)$, $\left(u_2,v_2\right)$, and $\left(u_3,v_3\right)$ being collinear:
\begin{Lemma}
  \label{lemma_determinant}
  If at least one of the elements $u_2,u_3,v_2,v_3\in\mathbb{Z}_{p^r}$ is invertible and Equation
  (\ref{eqn_collinear}) is fulfilled then the three points
  $\left(\overline{0},\overline{0}\right)$, $\left(u_2,v_2\right)$, $\left(u_3,v_3\right)$
  in $\mathbb{Z}_{p^r}^2$ are collinear.
\end{Lemma}
\begin{Proof}
  Due to symmetry we can assume that $u_2$ is invertible. Setting $t_1=\overline{1}$, $w_2=u_2$, $t_2=v_2u_2^{-1}$,
  and $w_3=u_3$, using the notation of Lemma \ref{lemma:fix_zero}, we obtain $v_3=u_3v_2u_2^{-1}$. Since this equation
  is equivalent to Equation (\ref{eqn_collinear}) we can use Lemma \ref{lemma:fix_zero} to conclude that
  $\left(\overline{0},\overline{0}\right)$, $\left(u_2,v_2\right)$, and $\left(u_3,v_3\right)$ are collinear.
\end{Proof}

\bigskip

Using the previous lemmas we can design an efficient algorithm, with some subroutines, to decide whether $r\ge 3$ points $\left(u_i,v_i\right)\in\mathbb{Z}_n^2$ are collinear or not. We assume that the prime factorization
$$
  n=\prod_{i=1}^l p_i^{\mu_i}
$$
is known in advance. For practical purposes our algorithm deals with integers instead of residue classes.
\begin{Algorithm}{\textbf{is\_collinear$\mathbf{\Big(\widehat{u}_1,\widehat{v}_1,\dots,\widehat{u}_r,
  \widehat{v}_r,p_1,\mu_1,\dots,p_l,\mu_l\Big)}$}\\}
  \label{algo_i}
  $n=\prod\limits_{i=1}^l p_i^{\mu_i}$\\
  $u_2=\widehat{u}_2-\widehat{u}_1$, $v_2=\widehat{v}_2-\widehat{v}_1$\\
  \texttt{if} $u_2\le 0$ \texttt{then} $u_2=u_2+n$\\
  \texttt{if} $v_2\le 0$ \texttt{then} $v_2=u_2+n$\\
  \texttt{for } $i$ \texttt{ from } $3$ \texttt{ to } $r$\\
  \hspace*{6mm}$u_3=\widehat{u}_i-\widehat{u}_1$, $v_2=\widehat{v}_i-\widehat{v}_1$\\
  \hspace*{6mm}\texttt{if} $u_3\le 0$ \texttt{then} $u_3=u_3+n$\\
  \hspace*{6mm}\texttt{if} $v_3\le 0$ \texttt{then} $v_3=u_3+n$\\
  \hspace*{6mm}\texttt{if } \textbf{is\_collinear\_fix\_zero$\mathbf{\left(u_2,v_2,u_3,v_3,
  p_1,\mu_1,\dots,p_l,\mu_l\right)}==false$}\\
  \hspace*{6mm}\texttt{then return} $false$\\
  \texttt{return} $true$
\end{Algorithm}

\begin{Algorithm}{\textbf{is\_collinear\_fix\_zero$\mathbf{\left(u_2,v_2,u_3,v_3,
  p_1,\mu_1,\dots,p_l,\mu_l\right)}$}\\}
  \label{algo_ii}
  \texttt{for } $i$ \texttt{ from } $1$ \texttt{ to } $l$\\
  \hspace*{6mm}$k=p_i^{\mu_i}$\\
  \hspace*{6mm}$u_2'=u_2-\left\lfloor\frac{u_2}{k}\right\rfloor\cdot k$\\
  \hspace*{6mm}$u_3'=u_3-\left\lfloor\frac{u_3}{k}\right\rfloor\cdot k$\\
  \hspace*{6mm}$v_2'=v_2-\left\lfloor\frac{v_2}{k}\right\rfloor\cdot k$\\
  \hspace*{6mm}$v_3'=v_3-\left\lfloor\frac{v_3}{k}\right\rfloor\cdot k$\\
  \hspace*{6mm}\texttt{if } \textbf{is\_collinear\_prime\_power$\mathbf{\left(u_2',v_2',u_3',v_3',
   p_i,\mu_i\right)}==false$}\\
  \hspace*{6mm}\texttt{then return} $false$\\
  \texttt{return} $true$
\end{Algorithm}

\begin{Algorithm}{\textbf{is\_collinear\_prime\_power$\mathbf{\left(u_2,v_2,u_3,v_3,
  p,r\right)}$}\\}
  \label{algo_iii}
  \texttt{if } $r==1$ \texttt{ then }\\
  \hspace*{6mm}\texttt{if } $u_2v_3\equiv u_3v_2\pmod{p}$ \texttt{ then return } $true$ \texttt{ else return } $false$\\
  \texttt{else }\\
  \hspace*{6mm}\texttt{if } $u_2\equiv u_3\equiv v_2\equiv v_3\equiv\overline{0}\pmod{p}$\\
  \hspace*{6mm}\texttt{then return } \textbf{is\_collinear\_prime\_power$\mathbf{\left(\frac{u_2}{p},\frac{v_2}{p},
  \frac{u_3}{p},\frac{v_3}{p},p,r-1\right)}$}\\
  \hspace*{6mm}\texttt{else }\\
  \hspace*{12mm}\texttt{if } $u_2v_3\equiv u_3v_2\pmod{p^r}$ \texttt{ then return } $true$ \texttt{ else return } $false$
\end{Algorithm}

Now we want to analyze the running time of Algorithm \ref{algo_i} and Algorithm \ref{algo_ii}:
\begin{Theorem}
  \label{thm:run_time}
  If $n=\prod\limits_{i=1}^lp_i^{\mu_i}$ then Algorithm \ref{algo_ii} needs at most
  $O\left(\sum\limits_{i=1}^l\mu_i\right)$
  time steps.
\end{Theorem}
\begin{Proof}
  It suffices to prove that Algorithm \ref{algo_iii} needs at most $r$ recursions, which is obvious.
\end{Proof}

\begin{Corollary}
  Given the prime factorization of $n$ Algorithm \ref{algo_ii} runs in $O(\log n)$ and Algorithm \ref{algo_i}
  runs in $O(r\cdot \log n)$ time.
\end{Corollary}

So let us have a (small) example to illustrate Algorithm \ref{algo_i}. We choose $n=625=5^4$, $u_1=\overline{1}$, $v_1=\overline{2}$, $u_2=\overline{76}$, $v_2=\overline{57}$, $u_3=\overline{251}$ and $v_3=\overline{102}$. At first we transform the problem to $u_1=v_1=\overline{0}$, $u_2=\overline{75}$, $v_2=\overline{55}$, $u_3=\overline{250}$, and $v_3=\overline{100}$. Since $n$ is a prime power we do not split up the problem into prime powers. Since the largest power of $5$ which divides all of $\widetilde{u}_2$, $\widetilde{u}_3$, $\widetilde{v}_2$, and $\widetilde{v}_3$ is $5^1$ we reduce the problem to $\left(\overline{0},\overline{0}\right)$, $\left(\overline{15},\overline{11}\right)$, $\left(\overline{50},\overline{20}\right)$ in $\mathbb{Z}_{125}$. Due to $\overline{15}\cdot\overline{20}=\overline{11}\cdot\overline{50}$ in $\mathbb{Z}_{125}$ the three original points are collinear.

If we have to check very often whether three points are collinear or not then it is more efficient to create a $\mathbb{Z}_n^2\times \mathbb{Z}_n^2$ table in order to bookmark whether $\left(\overline{0},\overline{0}\right)$, $p_1$, $p_2$ are collinear or not.

\section{Integer Linear Programming formulations}
\label{sec_ilp_formulations}

\noindent
In this section we formulate the problem of the exact determination of $m_2\!\left(\mathbb{Z}_n^2\right)$ as an integer linear program using the binary variables $x_{i,j}\in\{0,1\}$ for $1\le i,j\le n$. Here the variables $x_{i,j}$ encode a subset $$C:=\Big\{\left(i+\mathbb{Z}_n,j+\mathbb{Z}_n\right)\mid x_{i,j}=1,\,1\le i,j\le n\Big\}\subseteq\mathbb{Z}_n^2.$$
To enforce $C$ to be a cap, we require the linear inequality
\begin{equation}
  \label{eq:cap}
  \sum\limits_{i,j\,:\,\left(i+\mathbb{Z}_n,j+\mathbb{Z}_n\right)\in L} x_{i,j}\le 2
\end{equation}
for all lines $L$ of $\mathbb{Z}_n^2$. It is not difficult to show that these inequalities suffice to enforce that no three points of $C$ are collinear. Obviously we could also write up an inequality for every triple of collinear points, but we remark that Inequality (\ref{eq:cap}) is more compact. With the above we can state
\begin{alignat}{3}
  m_2\!\left(\mathbb{Z}_n^2\right)= &&\max\,\sum\limits_{i=1}^n\sum\limits_{j=1}^n x_{i,j}\label{ilp_max_cap}\\
  \text{subject to}\nonumber\\
  && \sum\limits_{i,j\,:\,\left(i+\mathbb{Z}_n,j+\mathbb{Z}_n\right)\in L} x_{i,j}\le 2\quad\quad
  \forall\text{ lines }L\text{ of }\mathbb{Z}_n^2\nonumber\\
  && x_{i,j}\in\{0,1\}\quad\quad\,\forall\, 1\le i,j\le n.\nonumber
\end{alignat}

\bigskip
\bigskip

We remark that every optimal solution of ILP (\ref{ilp_max_cap}) corresponds to a cap which is complete.

\noindent
Similarly we can state for the maximum cardinality of a cap which is a subset of a permutation point set:
\begin{alignat}{3}
  \sigma\!\left(\mathbb{Z}_n^2\right)= &&\max\,\sum\limits_{i=1}^n\sum\limits_{j=1}^n x_{i,j}\label{ilp_max_perm_cap}\\
  \text{subject to}\nonumber\\
  && \sum\limits_{i,j\,:\,\left(i+\mathbb{Z}_n,j+\mathbb{Z}_n\right)\in L} x_{i,j}\le 2\quad\quad
  \forall\text{ lines }L\text{ of }\mathbb{Z}_n^2\nonumber\\
  &&\sum\limits_{i=1}^n x_{i,j}\le 1\quad\,\,\quad\quad\forall 1\le j\le n\nonumber\\
  &&\sum\limits_{j=1}^n x_{i,j}\le 1\quad\,\,\quad\quad\forall 1\le i\le n\nonumber\\
  && x_{i,j}\in\{0,1\}\quad\quad\,\forall\, 1\le i,j\le n.\nonumber
\end{alignat}

\bigskip
\bigskip

\noindent
A bit more work is needed to express $n_2\!\left(\mathbb{Z}_n^2\right)$ as the optimal objective value of an ILP. Simply replacing $\max$ by $\min$ in ILP (\ref{ilp_max_cap}) would yield the optimal solution $x_{i,j}=0$ for all $1\le i,j\le n$. So we have to augment ILP (\ref{ilp_max_cap}) by some additional conditions and variables in order to enforce the, to the $x_{i,j}$ corresponding, set $C$ to be complete. Therefore we introduce line variables $y_L\in\{0,1\}$ for every line $L$ in $\mathbb{Z}_n^2$. The idea is that $y_L$ should equal $1$ if $C$ contains exactly two points of $L$. This can be modeled using the linear inequality
\begin{equation}
  \label{eq_y}
  1+y_L\ge\sum\limits_{i,j\,:\,\left(i+\mathbb{Z}_n,j+\mathbb{Z}_n\right)\in L} x_{i,j}\ge 2y_L
\end{equation}
for all lines $L$ in $\mathbb{Z}_n^2$.

To model the completeness of $C$ we introduce the linear inequality
\begin{equation}
  \label{eq_no_further_point}
  x_{i,j}+\sum\limits_{L\,:\,\left(i+\mathbb{Z}_n,j+\mathbb{Z}_n\right)\in L} y_L\ge 1
\end{equation}
for all $1\le i,j\le n$. Here the idea is, that a cap $C$ is complete if and only if there does not exist a point $P\in\mathbb{Z}_n^2\backslash C$ such that $C\cup\{P\}$ is also a cap. So let us assume that we have a binary variable allocation $x_{i,j}$, $y_L$ satisfying Inequalities (\ref{eq:cap}), (\ref{eq_y}), and (\ref{eq_no_further_point}), then for all $1\le i,j\le n$ either $P=\left(i+\mathbb{Z}_n,j+\mathbb{Z}_n\right)$ is contained in $C$ or there exist two points $P_1$, $P_2$ in $C$ such that $P_1$, $P_2$, and $P$ are collinear. Thus every feasible solution of the ILP
\begin{alignat}{3}
  n_2\!\left(\mathbb{Z}_n^2\right)= &&\min\,\sum\limits_{i=1}^n\sum\limits_{j=1}^n x_{i,j}\label{ilp_min_cap}\\
  \text{subject to}\nonumber\\
  && \sum\limits_{i,j\,:\,\left(i+\mathbb{Z}_n,j+\mathbb{Z}_n\right)\in L} x_{i,j}\le 2\quad\quad
  \forall\text{ lines }L\text{ of }\mathbb{Z}_n^2\nonumber\\
  && \sum\limits_{i,j\,:\,\left(i+\mathbb{Z}_n,j+\mathbb{Z}_n\right)\in L} x_{i,j}-2y_L\ge 0\quad
  \quad\forall\text{ lines }L\text{ of }\mathbb{Z}_n^2\nonumber\\
  && x_{i,j}+\sum\limits_{L\,:\,\left(i+\mathbb{Z}_n,j+\mathbb{Z}_n\right)\in L} y_L\ge 1\quad
  \quad\,\,\forall\, 1\le i,j\le n\nonumber\\
  && x_{i,j}\in\{0,1\}\quad\,\,\,\,\,\,\,\,\forall\, 1\le i,j\le n\nonumber\\
  && y_L\in\{0,1\}\quad\,\,\,\,\,\forall\text{ lines }L\text{ of }\mathbb{Z}_n^2\nonumber
\end{alignat}
corresponds to a complete cap $C$. For a given complete cap $C$ we can extend the corresponding partial variable allocation by setting $y_L=1$ exactly if $C$ contains exactly two points of $C$. Since we minimize the number of points of $C$ the target value of ILP (\ref{ilp_min_cap}) equals $n_2\!\left(\mathbb{Z}_n^2\right)$.

\medskip

Now let us have a look at the number of variables and inequalities of the ILPs (\ref{ilp_max_cap}), (\ref{ilp_max_perm_cap}), and (\ref{ilp_min_cap}). A function $\psi:\mathbb{N}\rightarrow\mathbb{N}$ is called \textbf{multiplicative} if $\psi(1)=1$ and  $\psi(nm)=\psi(n)\cdot\psi(m)$ for all coprime integers $n$ and $m$.
\begin{Definition}
  Let $\psi:\mathbb{N}\rightarrow\mathbb{N}$ be the multiplicative arithmetic  theoretic function defined by
  $\psi\left(p^r\right)=(p+1)p^{r-1}$ for prime powers $p^r>1$.
\end{Definition}
We remark that since
$$
  \lim\limits_{k\rightarrow\infty} \ln{p_k}\prod\limits_{i=1}^k \frac{1}{1+\frac{1}{p_k}}=\frac{\pi^2}{6e^\gamma},
$$
and
$$
  \lim\limits_{k\rightarrow\infty}\left(\prod\limits_{i=1}^k p_i\right)^{\frac{1}{p_k}}=e
$$
where $p_k$ denotes the $k$th prime and $\gamma$ denotes the Euler-Mascheroni constant with an approximate value of $0.57721566$, and due to $\psi\!\left(n^2\right)=n^2\cdot\prod\limits_{p|n}1+\frac{1}{p}$ (for primes $p$), we have
$$
  n^2\le \psi\!\left(n^2\right)\le 1.0828\cdot n^2\ln \ln n
$$
for all sufficiently large $n$.

\begin{Lemma}
  There are $\psi\!\left(n^2\right)$ lines in $\mathbb{Z}_n^2$ .
\end{Lemma}
\begin{Proof}
  Due to the Chinese remainder theorem the number of lines is a multiplicative arithmetic function. The number of cyclic
  subgroups of order $p^r$ in $\mathbb{Z}_{p^r}^2$ is given by
  $2\varphi\!\left(p^r\!\right)\left(p^r-\varphi\!\left(p^r\!\right)\right)+\varphi\!\left(p^r\!\right)^2=p^{2r}-p^{2r-2}$, where
  $\varphi\!\left(p^r\!\right)=p^{r-1}(p-1)$ is Euler's totient function. Since every cyclic subgroups of order $p^r$ contains
  $\varphi\!\left(p^r\!\right)=p^{r-1}(p-1)$ generators, there are $\frac{p^{2r}-p^{2r-2}}{p^{r-1}(p-1)}=(p+1)p^{r-1}$ lines through
  each point in $\mathbb{Z}_{p^r}^2$. In total there are $\frac{(p+1)p^{r-1}\cdot p^{2r}}{p^r}=\psi\!\left(p^{2r}\!\right)$
  lines since every line contains $p^r$ points.
\end{Proof}

Thus ILP (\ref{ilp_max_cap}) consists of $n^2$ variables and $\psi\!\left(n^2\right)$ inequalities, ILP (\ref{ilp_max_perm_cap}) consists of $n^2$ variables and $\psi\!\left(n^2\right)+2n$ inequalities, and ILP (\ref{ilp_min_cap}) consists of $n^2+\psi\!\left(n^2\right)$ variables and $2\psi\!\left(n^2\right)+n^2$ inequalities. So in all cases the number of variables and inequalities are in $O(n^2\ln n)$. But since generally the optimization (or also the feasibility problem) of $0$-$1$ linear programs is NP-complete these ILP formulations might not help too much from the theoretical point of view. On the other hand these ILP formulations enable us to determine some exact numbers and bounds of $m_2\!\left(\mathbb{Z}_n^2\right)$, $n_2\!\left(\mathbb{Z}_n^2\right)$, and $\sigma\!\left(\mathbb{Z}_n^2\right)$ in the next section.

\bigskip

\begin{figure}[htp]
  \begin{center}
    \setlength{\unitlength}{1.00cm}
    \begin{picture}(3,2)
      \put(0,0){\circle*{0.2}}
      \put(1,0){\circle*{0.2}}
      \put(2,0){\circle*{0.2}}
      \put(3,0){\circle*{0.2}}
      \put(0,1){\circle*{0.2}}
      \put(1,1){\circle*{0.2}}
      \put(2,1){\circle*{0.2}}
      \put(0,2){\circle*{0.2}}
      \put(1,2){\circle*{0.2}}
      \put(0,0){\line(1,0){3}}
      \put(0,0){\line(0,1){2}}
      \put(3,0){\line(-3,5){1.2}}
      \put(0,2){\line(1,0){1.8}}
    \end{picture}
    \quad\quad\quad\quad\quad\quad
    \begin{picture}(3,2)
      \put(0,0){\circle*{0.2}}
      \put(1,0){\circle*{0.2}}
      \put(2,0){\circle*{0.2}}
      \put(3,0){\circle*{0.2}}
      \put(0,1){\circle*{0.2}}
      \put(1,1){\circle*{0.2}}
      \put(2,1){\circle*{0.2}}
      \put(0,2){\circle*{0.2}}
      \put(1,2){\circle*{0.2}}
      \put(0,0){\line(1,0){3}}
      \put(0,0){\line(0,1){2}}
      \put(3,0){\line(-1,1){2}}
      \put(0,2){\line(1,0){1}}
    \end{picture}
  \end{center}
  \caption{Feasible sets of linear programs.}
  \label{fig:example_lp}
\end{figure}
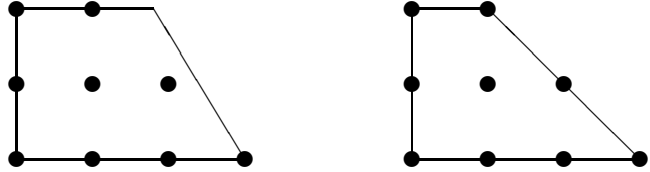

In contrast to ILP problems LP problems, i.~e.{} ILP problems without integrality constraints, can be solved in polynomial time. So in order to obtain an LP for the ILPs (\ref{ilp_max_cap}), (\ref{ilp_max_perm_cap}), and (\ref{ilp_min_cap}) we can relax the conditions $x_{i,j}\in\{0,1\}$, $y_L\in\{0,1\}$ by $0\le x_{i,j}\le 1$, $0\le y_L\le 1$. To solve the original (integral) problem several techniques, e.~g.{} branch~\&~bound, have to be applied. In many cases additional inequalities will be very useful for an optimization algorithm. We will explain this idea by considering the example
\begin{alignat}{3}
  &&\max\,x_1+2x_2\nonumber\\
  \text{subject to}\nonumber\\
  && 5x_1+3x_2\le 15\nonumber\\
  && x_2\le 2\nonumber\\
  && x_1,x_2\in\mathbb{N}_0\nonumber
\end{alignat}
On the left hand side of Figure \ref{fig:example_lp} we have depicted the feasible set of the relaxed linear program (i.~e.{} we have replaced $x_1,x_2\in\mathbb{N}_0$ by $x_1,x_2\ge 0$). The integral points are marked by filled circles. If we additionally require $x_1+x_2\le 3$ we obtain the feasible set as depicted on the right hand side of Figure \ref{fig:example_lp}.

We observe that the feasible set on the right hand side contains the same integral points as the feasible set of the left hand side, whereas the surface area differs. In this case we say that $x_1+x_2\le 3$ is a valid inequality w.~r.~t.{} the integral points. If the right hand side of Inequality (\ref{eq:cap}) in ILP (\ref{ilp_max_cap}) would be $1$ instead of $2$, then several classes of valid inequalities of the so called stable set polytope are known, e.~g.{}. odd circuit inequalities or clique inequalities (if given by edge inequalities).

Unfortunately we are not aware of any general (masking the geometric properties of $\mathbb{Z}_n^2$) results on valid inequalities for the polytope of ILP (\ref{ilp_max_cap}). 

\section{Bounds and exact values}
\label{sec_bounds_exact_values}

\noindent
In this section we want do state bounds and exact values for the three problems stated in the introduction which are either known or derived using an ILP based approach following the ILP formulations of Section \ref{sec_ilp_formulations}.

\subsection{Maximum cardinality of caps over $\mathbf{\mathbb{Z}_n^2}$}

\noindent
As mentioned in the introduction 
we have:
\begin{Lemma}
  $
    m_2\!\left(\mathbb{Z}_{nm}^2\right)\le\min\left\{n\cdot m_2\!\left(\mathbb{Z}_m^2\right),
  m_2\!\left(\mathbb{Z}_n^2\right)\cdot m\right\}
  $
  for coprime integers $n,m>1$.
\end{Lemma}
\begin{Proof}
  Let $X\subseteq\mathbb{Z}_{nm}^2$ be a cap, $l$ be a line in $\mathbb{Z}_n^2$, and
  $q_1:\mathbb{Z}_{nm}^2\rightarrow\mathbb{Z}_n^2$, $q_2:\mathbb{Z}_{nm}^2\rightarrow\mathbb{Z}_m^2$ be reduction maps. We set
  $Y:=\left\{x\in X\mid q_1(x)\in l\right\}$. If $|Y|\ge 2$ then $q_1(Y)\subseteq q_1(X)$ is collinear in $\mathbb{Z}_n^2$.
  Thus we must have $|Y|\le m_2\!\left(\mathbb{Z}_m^2\right)$ since otherwise $q_2(Y)$ is collinear in $\mathbb{Z}_m^2$, from which
  we could conclude the collinearity of $Y$ in $\mathbb{Z}_{nm}^2$. Since we can partition $\mathbb{Z}_n^2$ into $n$ lines
  and $m_2\!\left(\mathbb{Z}_m^2\right)\ge 2$ we have $|X|\le n\cdot m_2\!\left(\mathbb{Z}_m^2\right)$. Due to symmetry we 
  also have $|X|\le m\cdot m_2\!\left(\mathbb{Z}_n^2\right)$.
\end{Proof}

Directly solving ILP (\ref{ilp_max_cap}) using the commercial solver \texttt{ILOG CPLEX 11.2} yields the results given in Table \ref{table_direct}. We would like to remark that for $n\le 17$ the exact values of $m_2\!\left(\mathbb{Z}_n^2\right)$ were also determined in \cite{dipl_michael}.

\begin{table}[htp]
  \begin{center}
    \begin{tabular}{r|rrrrrrrrrrr}
      $n$ & 2 & 3 & 4 & 5 & 6 & 7 & 8 & 9 & 10 & 11 & 12 \\
      \hline $m_2\!\left(\mathbb{Z}_n^2\right)$ & 4 & 4 & 6 & 6 & 8 & 8 & 8 & 9 & 12 & 12 & 12\\
      \hline time in s & $<0.1$ & $<0.1$ & $<0.1$ & $<0.1$ & $<0.1$ & $0.65$ & $0.12$ & $48.05$ & $5.30$ &
      $126.24$ & $56.80$ \\
    \end{tabular}
    \caption{Directly solving ILP (\ref{ilp_max_cap}).}
    \label{table_direct}
  \end{center}
\end{table}

One of the main reasons why ILP solvers fail to tackle instances of ILP (\ref{ilp_max_cap}) for \textit{larger} values of $n$ is the highly symmetric formulation of the problem, i.~e.{} instead of few optimal solutions there exist a whole bulk of solutions which correspond to geometrically isomorphic caps. So we have to break symmetries by introducing further inequalities. To this end we consider automorphisms $\alpha$ of $\mathbb{Z}_n^2$, i.~e.{} mappings from $\mathbb{Z}_n^2$ to $\mathbb{Z}_n^2$ which preserve 
point-line incidences. This especially means that $C\subseteq\mathbb{Z}_n^2$ is a cap if and only if $\alpha(C)$ 
is a cap. Since the translations $x\mapsto x+y$ for $y\in\mathbb{Z}_n^2$ are such automorphisms we can assume w.l.o.g.{} that the point $\left(\overline{0},\overline{0}\right)$ is contained in each non-empty cap $C$. Another class of automorphisms arises from multiplication of invertible $2\times 2$~matrices of $\mathbb{Z}_n$. So let us denote the resulting group of all translations and all invertible $2\times 2$~matrices by $G$. (All elements of $G$ are automorphisms.)

In the following we will assume that some points $a_i\in\mathbb{Z}_n^2$ are not contained in cap $C$ but some points $b_i\in\mathbb{Z}_n^2$ are contained in $C$. From this we can conclude some valid inequalities for our ILPs using the group $G$. If there exists an index $h$, an element $z\in\mathbb{Z}_n^2$, and an automorphism $\alpha$ such that
$$
  \alpha\left(\left\{z,a_1,\dots\right\}\right)=\left\{b_h,a_1,\dots\right\}
$$
then we can conclude w.l.o.g.{} that $z\notin C$. In the formulation as an integer linear program this translates to an equation $x_{i,j}=0$ (for $\left(i+\mathbb{Z}n,j+\mathbb{Z}n\right)=z$). If there exists an index $h$, three elements $z_1,z_2,z_3\in\mathbb{Z}_n^2$, and an automorphism $\alpha$ such that
$$
  \alpha\left(\left\{z_1,z_2,a_1,\dots\right\}\right)=\left\{b_h,z_3,a_1,\dots\right\}
$$
then we can conclude w.l.o.g.{} that not both of $z_1$ and $z_2$ are contained in $C$. In the ILP formulation this translates to an inequality $x_{i,j}+x_{i',j'}\le 1$ (for $\left(i+\mathbb{Z}n,j+\mathbb{Z}n\right)=z_1$ and $\left(i'+\mathbb{Z}n,j'+\mathbb{Z}n\right)=z_2$).

With this tool at hand we are able to determine some further exact values of $m_2\!\left(\mathbb{Z}_n^2\right)$ via case differentiations and the \texttt{ILOG CPLEX} solver, see Table \ref{table_case_dif}.

\begin{table}[htp]
  \begin{center}
    \begin{tabular}{r|l}
      $n$ & cases\\
      \hline
      14 & $a_1=\left(\overline{0},\overline{0}\right)$, $b_1=\left(\overline{1},\overline{0}\right)$
                        resolved in $0.57$~s.\\
         & $a_1=\left(\overline{0},\overline{0}\right)$, $a_2=\left(\overline{1},\overline{0}\right)$,
                        $b_1=\left(\overline{0},\overline{1}\right)$ resolved in $3.21$~s.\\
         & $a_1=\left(\overline{0},\overline{0}\right)$, $a_2=\left(\overline{1},\overline{0}\right)$,
                        $a_3=\left(\overline{0},\overline{1}\right)$ resolved in $675.54$~s.\\
      \hline
      15 &  $a_1=\left(\overline{0},\overline{0}\right)$, $b_1=\left(\overline{1},\overline{0}\right)$
                        resolved in $0.04$~s. \\
         & $a_1=\left(\overline{0},\overline{0}\right)$, $a_2=\left(\overline{1},\overline{0}\right)$,
                        $b_1=\left(\overline{0},\overline{1}\right)$ resolved in $0.67$~s.\\
         & $a_1=\left(\overline{0},\overline{0}\right)$, $a_2=\left(\overline{1},\overline{0}\right)$,
                        $a_3=\left(\overline{0},\overline{1}\right)$ resolved in $1830.68$~s.\\
      \hline 
      16 &$a_1=\left(\overline{0},\overline{0}\right)$, $b_1=\left(\overline{1},\overline{0}\right)$
                        resolved in $0.17$~s.\\
          & $a_1=\left(\overline{0},\overline{0}\right)$, $a_2=\left(\overline{1},\overline{0}\right)$,
                        $b_1=\left(\overline{0},\overline{1}\right)$ resolved in $0.54$~s.\\
                  & $a_1=\left(\overline{0},\overline{0}\right)$, $a_2=\left(\overline{1},\overline{0}\right)$,
                        $a_3=\left(\overline{0},\overline{1}\right)$ resolved in $251.93$~s.\\
      \hline
      18 & $a_1=\left(\overline{0},\overline{0}\right)$, $b_1=\left(\overline{1},\overline{0}\right)$
                        resolved in $0.07$~s.\\
       & $a_1=\left(\overline{0},\overline{0}\right)$, $a_2=\left(\overline{1},\overline{0}\right)$,
                        $b_1=\left(\overline{0},\overline{1}\right)$ resolved in $0.13$~s.\\
             & $a_1=\left(\overline{0},\overline{0}\right)$, $a_2=\left(\overline{1},\overline{0}\right)$,
                        $a_3=\left(\overline{0},\overline{1}\right)$ resolved in $2554.91$~s.\\
     \hline 
     20 &  $a_1=\left(\overline{0},\overline{0}\right)$, $b_1=\left(\overline{1},\overline{0}\right)$
                        resolved in less than $1$~s.\\
                  & $a_1=\left(\overline{0},\overline{0}\right)$, $a_2=\left(\overline{1},\overline{0}\right)$,
                        $b_1=\left(\overline{0},\overline{1}\right)$ resolved in less than $1$~s.\\
                  & $a_1=\left(\overline{0},\overline{0}\right)$, $a_2=\left(\overline{1},\overline{0}\right)$,
                        $a_3=\left(\overline{0},\overline{1}\right)$, $a_4=\left(\overline{4},\overline{1}\right)$
                        resolved in $1950$~s.\\
                  & $a_1=\left(\overline{0},\overline{0}\right)$, $a_2=\left(\overline{1},\overline{0}\right)$,
                        $a_3=\left(\overline{0},\overline{1}\right)$, $a_4=\left(\overline{4},\overline{4}\right)$,
                        $b_1=\left(\overline{4},\overline{1}\right)$ resolved in $513$~s.\\
                  & $a_1=\left(\overline{0},\overline{0}\right)$, $a_2=\left(\overline{1},\overline{0}\right)$,
                        $a_3=\left(\overline{0},\overline{1}\right)$, $a_4=\left(\overline{5},\overline{4}\right)$,
                        $b_1=\left(\overline{4},\overline{1}\right)$,\\& $b_2=\left(\overline{4},\overline{4}\right)$
                        resolved in $80$~s.\\
                  & \dots, $a_4=\left(\overline{5},\overline{3}\right)$, \dots resolved in less than $1$~s.\\
                  & \dots, $a_4=\left(\overline{9},\overline{1}\right)$, \dots resolved in less than $1$~s.\\
                  & \dots, $a_4=\left(\overline{16},\overline{9}\right)$, \dots resolved in $3584$~s. \\
                  & \dots, $a_4=\left(\overline{15},\overline{11}\right)$, \dots resolved in $800$~s. \\
                  & \dots, $a_4=\left(\overline{15},\overline{10}\right)$, \dots resolved in $3562$~s.\\
                  & \dots, $a_4=\left(\overline{6},\overline{5}\right)$, \dots resolved in $3725$~s.\\
                  & \dots, $a_4=\left(\overline{5},\overline{5}\right)$, \dots resolved in $5922$~s.\\
                  & \dots, $a_4=\left(\overline{10},\overline{1}\right)$, \dots resolved in less than $1$~s.\\
                  & \dots, $a_4=\left(\overline{5},\overline{1}\right)$, \dots resolved in $3080$~s.\\
                  & $a_1=\left(\overline{0},\overline{0}\right)$, $a_2=\left(\overline{1},\overline{0}\right)$,
                        $a_3=\left(\overline{0},\overline{1}\right)$, $b_1=\left(\overline{4},\overline{1}\right)$,
                        $b_2=\left(\overline{4},\overline{4}\right)$, $b_3=\left(\overline{5},\overline{4}\right)$,\\&
                        $b_4=\left(\overline{5},\overline{3}\right)$, $b_5=\left(\overline{9},\overline{1}\right)$,
                        $b_6=\left(\overline{16},\overline{9}\right)$, $b_7=\left(\overline{15},\overline{11}\right)$,
                        $b_8=\left(\overline{15},\overline{10}\right)$, $b_9=\left(\overline{6},\overline{5}\right)$,\\&
                        $b_{10}=\left(\overline{5},\overline{5}\right)$,
                        $b_{11}=\left(\overline{10},\overline{1}\right)$,
                        $b_{12}=\left(\overline{5},\overline{1}\right)$ resolved in less than $1$~s.\\
    \hline
   \end{tabular}
    \caption{Case differentiations for the determination of $m_2\!\left(\mathbb{Z}_n^2\right)$.}
    \label{table_case_dif}
  \end{center}
\end{table}

In each case we consider ILP (\ref{ilp_max_cap}) augmented by the inequalities arising from the $a_i$, $b_i$ as described above and by the inequality $\sum\limits_{i=1}^n\sum\limits_{j=1}^n x_{i,j}\ge l+1$, where $l$ is the cardinality of a cap $C$ in $\mathbb{Z}_n^2$. Since it is not hard to find a cap of cardinality $m_2\!\left(\mathbb{Z}_n^2\right)$ (one can e.~g.{} use an ILP solver), we can choose $l=m_2\!\left(\mathbb{Z}_n^2\right)$. Thus \textit{resolved} means that the ILP solver has proven that no integer solution can exist in each of the stated subcases. In Table \ref{table_m_cases} we give the proven exact values of $m_2\!\left(\mathbb{Z}_n^2\right)$ and some bounds which can be obtained by applying the described methods.

\begin{table}[htp]
  \begin{center}
    \begin{tabular}{r|rrrrrrrr}
      $n$ & 14 & 15 & 16 & 18 & 20 & 21 & 22 & 24\\
      \hline $m_2\!\left(\mathbb{Z}_n^2\right)$ & 12 & 15 & 14 & 17 & 18 & 18 & 18--24 & 18--24\\
    \end{tabular}
    \caption{Solving ILP (\ref{ilp_max_cap}) utilizing case differentions and symmetry.}
    \label{table_m_cases}
  \end{center}
\end{table}

\begin{figure}[htp]
  \begin{center}
    \setlength{\unitlength}{0.20cm}
    \begin{picture}(25.2,25.2)
      \multiput(0.15,0)(0,1){26}{\line(1,0){25}}
      \multiput(0.15,0)(1,0){26}{\line(0,1){25}}
      \put(0.65,0.5){\circle*{0.6}}   
      \put(1.65,0.5){\circle*{0.6}}   
      \put(13.65,5.5){\circle*{0.6}}  
      \put(5.65,6.5){\circle*{0.6}}   
      \put(17.65,6.5){\circle*{0.6}}  
      \put(7.65,7.5){\circle*{0.6}}   
      \put(13.65,7.5){\circle*{0.6}}  
      \put(18.65,8.5){\circle*{0.6}}  
      \put(19.65,8.5){\circle*{0.6}}  
      \put(16.65,9.5){\circle*{0.6}}  
      \put(7.65,10.5){\circle*{0.6}}  
      \put(9.65,10.5){\circle*{0.6}}  
      \put(3.65,11.5){\circle*{0.6}}  
      \put(21.65,11.5){\circle*{0.6}} 
      \put(4.65,12.5){\circle*{0.6}}  
      \put(2.65,13.5){\circle*{0.6}}  
      \put(10.65,13.5){\circle*{0.6}} 
      \put(15.65,14.5){\circle*{0.6}} 
      \put(18.65,14.5){\circle*{0.6}} 
      \put(21.65,23.5){\circle*{0.6}} 
    \end{picture}
  \end{center}
  \caption{A cap of cardinality $20$ over $\mathbb{Z}_{25}^2$.}
  \label{fig_ex_25}
\end{figure}
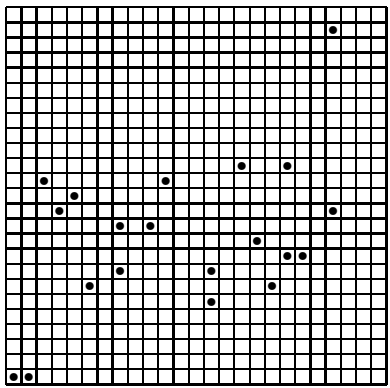

\noindent
The starting point of our studies was the determination of $m_2\!\left(\mathbb{Z}_{25}^2\right)$. Is is not too hard to find a cap of cardinality $20$ in $\mathbb{Z}_{25}^2$, see Figure \ref{fig_ex_25} for an example. In the projective case very recently, a bit surprising, a cap of cardinality $21$ was found \cite{matthias}. Soon after, Kohnert et al.{}, see \cite{hp_algorithm}, verified this constructive result by prescribing a cyclic group of order $3$ as a subgroup of the automorphism group of a cap. For the general method of prescribing automorphisms we refer e.~g.{} to \cite{gm2,gm1}. 

We remark that the known bounds for the maximum size of a cap in the projective Hjelmslev plane $\operatorname{PHG}\left(\mathbb{Z}_{25}^3\right)$ were $20\dots 25$, see e.~g.{} \cite{0994.51006}. Here we conjecture $21$ to be the correct value.
It is very likely that the method of M.~Koch can be continued to an exhaustive search to resolve this case.

%
%
%

To determine $m_2\!\left(\mathbb{Z}_{25}^2\right)=20$ using the ILP based approach one has to perform several case differentions, use some geometric insights and exhaustive enumeration. At first we look at pairs $a,b\in\mathbb{Z}_{25}^2$ of points which lie in the same neighborhood $\sim_1$, i.~e.{} where $a$ and $b$ are incident with $5^1$ common lines, see Section \ref{sec_collinear}. In $\mathbb{Z}_{25}^2$ two points $\left(a_1,a_2\right)$, $\left(b_1,b_2\right)$ belong to such a class if and only if both $\widehat{a}_1-\widehat{b}_1$ and $\widehat{a}_2-\widehat{b}_2$ are divisible by $5$. If three pairwise different points $c_1,c_2,c_3\in\mathbb{Z}_{25}^2$ fulfill $a\sim_1 b\sim_1 c$ then the caps $C$ containing $c_1$, $c_2$, and $c_3$ are rather small, i.~e.{} $|C|<20$, which can be easily shown both theoretically and computationally (by solving an ILP).\footnote{In the printed version of this paper, we have incorrectly asserted, that in this case the three points are collinear. Indeed, one can even show $|C|\le 14$ in this case.}
Now we consider caps $C$ which contain at least three pairs $c_1\sim_1 c_2$, $c_3\sim_1 c_4$, and $c_5\sim_1 c_6$. Since $C$ is a cap we have $c_1\not\sim_1 c_3$, $c_1\not\sim_1 c_5$, and $c_3\not\sim_1 c_5$. Due to the symmetry group $G$, defined above as a subgroup of the automorphism group, we can assume w.l.o.g.{} $c_1=\left(\overline{0},\overline{0}\right)$, $c_3=\left(\overline{1},\overline{0}\right)$, and $c_5=\left(\overline{0},\overline{1}\right)$. For $\{c_1,\dots,c_6\}^G$ there are $104$~orbits.
In each of these $104$~cases the augmented ILP (\ref{ilp_max_cap}) has an optimal solution less than $21$.

Thus we only have to consider caps with at most two pairs $c_1\sim_1 c_2$ and $c_3\sim_1 c_4$. So in the next step we consider all $18$~orbits 
$\{c_1,\dots,c_5\}^G$ with $c_1=\left(\overline{0},\overline{0}\right)$, $c_3=\left(\overline{1},\overline{0}\right)$, $c_1\sim_1 c_4$, and $c_2\sim_1 c_5$. As additional inequalities we can use that no further pair $c_i\sim_1 c_j$ can exist in $C$. Also, here the optimal solution of the augmented ILPs (\ref{ilp_max_cap}) is less than $21$.

In a further step we consider all $39$ orbits 
$\{c_1,\dots,c_4\}^G$ with $c_1=\left(\overline{0},\overline{0}\right)$, $c_3=\left(\overline{1},\overline{0}\right)$, and $c_1\sim_1 c_4$. Again the optimal solution of the augmented ILPs (\ref{ilp_max_cap}) is less than $21$. Thus in a possible cap of cardinality at least $21$ in $\mathbb{Z}_{25}^2$ no pair $c_i\sim_1 c_j$ can exist.

Next we look at the orbits of $\{c_1,\dots,c_4\}^G$ with $c_1=\left(\overline{0},\overline{0}\right)$, $c_3=\left(\overline{1},\overline{0}\right)$ where no pair $c_i\sim_1 c_j$ occurs. We pick a representative of each orbit and label them by $z_1,\dots,z_{33}$. 
For $i\ge 4$ we set $a_1=c_1$, $a_2=c_2$, $a_3=c_3$, $a_4=z_i$, $b_1=z_1$, \dots, $b_{i-1}=z_{i-1}$ in order to prescribe or to forbid some points of cap $C$. In each of these cases we can verify that $C$ must contain less than $21$ elements. Thus only three possibilities for $\{c_1,\dots,c_4\}^G$ are left.

For the final step we utilize orderly generation \cite{0392.05001} with some look ahead: we run through all $25$ equivalence classes $\left\{x\in\mathbb{Z}_{25}\mid x\sim_1 y\right\}$. In each equivalence class we can pick at most one point. So either we pick one element of such an equivalence class or we decide not to take an element of this equivalence class. If we have fixed some elements of our cap $C$ then it may happen that there are equivalence classes where we have no possibility to pick an element due to the restriction of at most $2$ points on a line. So for each partial cap we can count the number $r$ of equivalence classes, where it is possible to select a further point for cap $C$. If $|C|+r<21$ then we can stop extending cap $C$.




\subsection{Minimum cardinality of complete caps over $\mathbf{\mathbb{Z}_n^2}$}
From \cite{pre05050776} we can cite the bounds $$\max\left(4,\sqrt{2p}+\frac{1}{2}\right)\le n_2\!\left(\mathbb{Z}_{n}^2\right)\le\max(4,p+1),$$ where $p$ is the smallest prime divisor of $n$, and
$$n_2\!\left(\mathbb{Z}_{nm}^2\right)\le \min\Big(n_2\!\left(\mathbb{Z}_{n}^2\right),n_2\!\left(\mathbb{Z}_{m}^2\right)\Big)$$ for coprime integers $n$ and $m$. Additionally the author of \cite{pre05050776} conjectures $n_2\!\left(\mathbb{Z}_{p^a}^2\right)\le n_2\!\left(\mathbb{Z}_{p^b}^2\right)$ for $a\le b$, $p$ being a prime and $n_2\!\left(\mathbb{Z}_{nm}^2\right)=\min\Big(n_2\!\left(\mathbb{Z}_{n}^2\right),n_2\!\left(\mathbb{Z}_{m}^2\right)\Big)$ for coprime integers $n$ and $m$. If these conjectures turn out to be true then it would suffice to determine the values $n_2\!\left(\mathbb{Z}_{p}^2\right)$ for primes $p$.

If $n$ is divisible by $2$ or $3$ then we can conclude $n_2\!\left(\mathbb{Z}_{n}^2\right)=4$ from the above inequalities and $n_2\!\left(\mathbb{Z}_{2}^2\right)=n_2\!\left(\mathbb{Z}_{3}^2\right)=4$. For $n>1$ we have $n_2\!\left(\mathbb{Z}_{n}^2\right)\ge 4$. So in Table \ref{table_n_2} and Table \ref{table_n_2_2} we have given the exact values and bounds for $n_2\!\left(\mathbb{Z}_{n}^2\right)$ arising from ILP (\ref{ilp_max_cap}), where $n$ is either a prime or coprime to $6$. We would like to remark that some of these numbers are already given in \cite{pre05050776}. If $n$ is a prime then we can assume w.l.o.g.{} that the points $\left(\overline{0},\overline{0}\right)$, $\left(\overline{1},\overline{0}\right)$, and $\left(\overline{0},\overline{1}\right)$ are contained in cap $C$. In this case some results from \cite{1115.51005} can be used, e.~g.{} we have $n_2\!\left(\mathbb{Z}_{n}^2\right)\ge t(2,q)-2$,
where $t(2,q)$ denotes the smallest size of a complete cap in $PG(2,q)$. For the exact values of $t(2,q)$ for $q\le 29$ we refer to \cite{1025.51012,1044.51007}. In \cite{1115.51005} there is mentioned a construction that produces a cap of size $\frac{q-3}{2}$ which is complete at least for $q>413$. The exact value of $n_2\!\left(\mathbb{Z}_{25}^2\right)$ is known to be $6$ but unfortunately
we were not able to validate the lower bound using our ILP based approach.

\begin{table}[htp]
  \begin{center}
    \begin{tabular}{r||c|c|c|c|c|c|c|c|c|c|c}
      $n$ & 2 & 3 & 5 & 7 & 11 & 13 & 17 & 19 & 23 & 25 & 29\\
      \hline
      $n_2\!\left(\mathbb{Z}_{n}^2\right)$ &
      $4$ & $4$ & $5$ & $6$ & $7$ & $8$ & $8\dots10$ & $8\dots11$ & $8\dots12$ & $4\dots 6$ & $11\dots15$ \\
    \end{tabular}
    \caption{Values of $n_2\!\left(\mathbb{Z}_{n}^2\right)$ for small $n$ which are either prime or coprime to $6$.}
    \label{table_n_2}
  \end{center}
\end{table}

\begin{table}[htp]
  \begin{center}
    \begin{tabular}{r||c|c|c|c|c}
      $n$ & 31 & 37 & 41 & 43 & 47\\
      \hline
      $n_2\!\left(\mathbb{Z}_{n}^2\right)$ &
      $9\dots16$ & $10\dots17$ & $10\dots20$ & $10\dots21$ & $11\dots22$\\
    \end{tabular}
    \caption{Values of $n_2\!\left(\mathbb{Z}_{n}^2\right)$ for small $n$ which are either prime or coprime to $6$.}
    \label{table_n_2_2}
  \end{center}
\end{table}



\subsection{Maximum cardinality of caps over $\mathbf{\mathbb{Z}_n^2}$ which are subsets of permutations}
Obviously we have $\sigma\!\left(\mathbb{Z}_{n}^2\right)\le n$ since a permutation consists of $n$ points. Utilizing the ILP formulation (\ref{ilp_max_perm_cap}) and the \texttt{ILOG CPLEX} solver we have obtained the values and bounds of Table \ref{table_perm} and Table \ref{table_perm2}. So e.~g.{} for $n\in\{1,2,4,6,8,12\}$ there exist permutations whose graphs in $\mathbb{Z}_n^2$ are caps. We would like to remark that for this problem the applicable group of automorphisms is much smaller then for the other two problems. Here translations, changes of the coordinate axes, and reflection at one of the axes are automorphisms.

\begin{table}[htp]
\begin{center}
\begin{tabular}{c|ccccccccccccccccccc}
 $n$ & 1 & 2 & 3 & 4 & 5 & 6 & 7 & 8 & 9 & 10 & 11 & 12 & 13 & 14 & 15 & 16 & 17 & 18 & 19\\
 \hline
 $\sigma(n)$&1&2&2&4&4&6&6&8&6&8&10&12&12&12&13&13&16&13&18\\
\end{tabular}
\caption{Values and bounds for $\sigma\!\left(\mathbb{Z}_{n}^2\right)$.}
\label{table_perm}
\end{center}
\end{table}

\begin{table}[htp]
\begin{center}
\begin{tabular}{c|ccccccccccc}
 $n$ & 20 & 21 & 22 & 23 & 24 & 25 & 26 & 27 & 28 & 29 & 30\\
 \hline
 $\sigma(n)$&16&16-17&16-17&22&20-22&19-22&18-24&18-25&22-27&28&22-29\\
\end{tabular}
\caption{Values and bounds for $\sigma\!\left(\mathbb{Z}_{n}^2\right)$.}
\label{table_perm2}
\end{center}
\end{table}

\end{document}